\documentclass[aos,MSNbibl,seceqn,rotating,dvips]{arximspdf}
\usepackage{dcolumn}

\doi{10.1214/13-AOS1180} 
\volume{42}
\issue{1}
\pubyear{2014}
\firstpage{64}
\lastpage{86}

\makeatletter
\newcolumntype{d}[1]{D{.}{.}{#1}}
\newcommand{\rrvert}{\vert}
\newcommand{\llvert}{\vert}
\def\cal{\mathcal}

\def\convinlaw{\stackrel{{\cal L}}{\Longrightarrow}}
\def\convinp{\stackrel{P}{\longrightarrow}}
\def\RR{\mathbb{R}}
\def\ZZ{\mathbb{Z}}
\def\EE{\mathbb{E}}
\def\WW{\mathbb{W}}
\def\YY{\mathbb{Y}}
\def\tends{\rightarrow}
\newproclaim{ass}{Assumptions}

\newtheorem{thmm}{Theorem}[section]
\newtheorem{Lemma}{Lemma}[section]

\newtheorem{Proposition}{Proposition}[section]
\newproclaim{Definition}{Definition}[section]
\newproclaim{Remark}{Remark}[section]
\makeatother

\begin{document}
\begin{frontmatter}

\title{Asymptotic theory of cepstral random fields}
\runtitle{Cepstral random fields}

\begin{aug}
\author[a]{\fnms{Tucker S.}~\snm{McElroy}\ead[label=e1]{tucker.s.mcelroy@census.gov}}
\and
\author[b]{\fnms{Scott~H.} \snm{Holan}\corref{}\thanksref{t3}\ead[label=e2]{holans@missouri.edu}}

\thankstext{t3}{Supported by NSF and the U.S. Census Bureau under NSF Grant
SES-1132031, through the NSF-Census Research Network (NCRN)
program and by a University of Missouri Research Board grant.}
\runauthor{T.~S. McElroy and S.~H. Holan}

\affiliation{U.S. Census Bureau and University of
Missouri}

\address[a]{Center for Statistical Research and Methodology\\
U.S. Census Bureau\\
4600 Silver Hill Road\\
Washington, D.C. 20233-9100\\
USA\\
\printead{e1}}

\address[b]{Department of Statistics\\
University of Missouri\\
146 Middlebush Hall\\
Columbia, Missouri 65211-6100\\
USA\\
\printead{e2}}
\end{aug}

\received{\smonth{3} \syear{2013}}
\revised{\smonth{10} \syear{2013}}

%
\begin{abstract}
Random fields play a central role in the analysis of spatially
correlated data and, as a result, have a significant impact on a
broad array of scientific applications. This paper studies the
cepstral random field model, providing recursive formulas that
connect the spatial cepstral coefficients to an equivalent
moving-average random field, which facilitates easy computation of
the autocovariance matrix. We also provide a comprehensive
treatment of the asymptotic theory for two-dimensional random field
models: we establish asymptotic results for Bayesian, maximum
likelihood and quasi-maximum likelihood estimation of random field
parameters and regression parameters. The theoretical results are
presented generally and are of independent interest, pertaining to a
wide class of random field models. The results for the cepstral
model facilitate model-building: because the cepstral coefficients
are unconstrained in practice, numerical optimization is greatly
simplified, and we are always guaranteed a positive definite
covariance matrix. We show that inference for individual
coefficients is possible, and one can refine models in a disciplined
manner. Our results are illustrated through simulation and the
analysis of straw yield data in an agricultural field experiment.
\end{abstract}

%

\end{frontmatter}

\section{Introduction}\label{sec:intro}
Spatial data feature heavily in many scientific disciplines
including ecology, environmental science, epidemiology, geography,
geology, small area estimation, and socio-demographics. Although
spatial data can be broadly placed into three categories:
\textit{geostatistical data}, \textit{lattice data} and
\textit{spatial patterns}~\cite{cressie1993statistics}, our focus
mainly resides in the development of cepstral random field models
for spatial lattice data. That is, we consider random fields where
the index set for the variables is $\ZZ^2$, appropriate for image
processing, for example.

Research on spatial random fields dates back over half a century; for
example, see Whittle \cite{whittle1954}. Other references on spatial random
fields include Besag \cite{besag1972correlation,besag1974spatial},
Guyon \cite{guyon1982parameter}, Rosenblatt \cite{murray1985stationary},
Besag and Green \cite{besag1993spatial} and Rosenblatt \cite
{rosenblatt2000gaussian}, among
others. Comprehensive overviews can be found in
Cressie \cite{cressie1993statistics}, Stein \cite{stein1999interpolation},
Banerjee, Carlin and
Gelfand \cite{banerjee2004hierarchical}, Cressie and Wikle \cite
{cressiewikle2011} and the
references therein. Recently, there has been a growing interest in
modeling spatial random fields through the spectral domain. For
example, see Fuentes \cite{fuentes2002spectral}, Fuentes, Guttorp and
Sampson~\cite{fuentes2007using},
Tonellato \cite{tonellato2007random}, Fuentes and Reich \cite
{fuentes2010spectral}, Bandyopadhyay and
Lahiri \cite{bandyopadhyayasymptotic} and the
references therein.

For a stationary Gaussian random field, it is natural to impose a
Markov structure, as described in Rue and Held \cite{rueheld:2010}, in
order to
obtain an inverse covariance matrix (i.e., a precision matrix) that
has a sparse structure, because this will ensure speedy computation
of maximum likelihood estimates. Rue and Held \cite{rueheld:2010} show how
careful specification of conditional distributions generates a
well-defined random field. However, this technique relies upon
imposing a priori a sparse structure on the precision matrix,
that is, demanding that many conditional precisions be zero. In
contrast, the cepstral random field does not generate a sparse
covariance (or precision) matrix, and yet always yields a
well-defined spatial random field; this occurs because the model is
formulated in the frequency domain by ensuring a positive spectral
density. This frequency-domain approach provides a general way of
specifying a nonisotropic random field, which is useful when we do
not have a~prior notion about conditional variances or precisions.

The cepstral random field allows for unconstrained optimization of
the objective function, that is, each model coefficient can be any real
number independently of the others; this appealing property is in
marked contrast to other models and approaches, such as moving
averages or Markov random fields (these require constraints on
parameters to achieve identifiability and/or a well-defined
process).
In the development of this model, Solo \cite{solo1986modeling} presents
estimation approaches by both log periodogram regression and Whittle
maximum likelihood, but does not derive the asymptotic properties of
estimators. Based on information criterion, Mallows's $C_p$, and
hypothesis testing, the author briefly describes methods for model
selection. Some key advantages of the cepstral model are that it is
well defined (because it is defined through the spectral density),
it is identifiable and the cepstral parameter estimates are
asymptotically uncorrelated with one another.

This paper provides a
first comprehensive treatment of the theory for cepstral random
field models. In particular, we establish recursive formulas for
connecting cepstral random fields to moving average random fields,
thus facilitating efficient computation of the spatial
autocovariances, which are needed for likelihood evaluation and
prediction. Critically, the resulting autocovariance matrix is
guaranteed to be postive-definite; note that if we were to work with
a moving average (MA) field instead, it would not be identifiable
without imposing further complicated parameter restrictions.

Additionally, we develop asymptotic results for Bayesian, maximum
likelihood, and quasi-maximum likelihood estimation of field
parameters and regression parameters under an expanding domain
formulation. In particular, we establish asymptotic consistency in
both the Bayesian and likelihood settings and provide central limit
theorems for the frequentist estimators we propose. We discuss the
computational advantages of the cepstral model, and propose an
exact Whittle likelihood that avoids the burdensome inversion of the
autocovariance matrix. Although our primary focus is on cepstral
models, the theoretical developments are presented for general
random field models with regression effects. Our results are of
independent interest and extend the existing results of
Mardia and Marshall~\cite{mardia1984maximum}, providing a rigorous
framework for
conducting model building and inference under an expanding domain
framework; this is applicable to lattice random field data that is
sampled at regular fixed intervals, and for which in-filling is
either impractical or of little interest.

As discussed in Sections~\ref{sec:model} and \ref{sec:fitting}, the
proposed cepstral models are computationally advantageous over many
current models (e.g., spatial autoregressive models), because no
constraints need to be imposed on the parameters to ensure the
resulting autocovariance matrix remains positive definite. In fact,
given the recursive formulas of Section~\ref{sec:model}, one can
model the two-dimensional cepstral coefficients (i.e., the Fourier
coefficients of the two-dimensional log spectrum) and arrive at the
autocovariances without the need for direct Fourier inversion.

Since the model's first inception \cite{solo1986modeling}, the
cepstral random field literature has remained sparse, with
relatively few examples to date. For example, Cressie~\cite{cressie1993statistics},
page~448, makes brief mention of the model. In a
different context, Noh and Solo~\cite{noh2007true} use cepstral random
fields to
test for space--time separability. Sandgren and
Stoica \cite{sandgren2006nonparametric}
use two-dimensional cepstrum thresholding models to estimate the
two-dimensional spectral density. However, this work does not treat
the random field case. Related to our work,
Kizilkaya and
Kayran \cite{kizilkaya2005arma} derive an algorithm for computing cepstral
coefficients from a known ARMA random field, whereas
Kizilkaya \cite{kizilkaya2007parameter} provides a recursive formula for
obtaining nonsymmetric half plane MA random field models for a
given cepstral specification. In contrast, our recursive formulas
provide unrestricted MA random fields as well as the necessary
autocovariances for expressing the Gaussian likelihood.

This paper proceeds as follows. Section~\ref{sec:model} describes
the cepstral model and its computation. Specifically, this section
lays out the recursive formulas that are needed to estimate the
autocovariances given the cepstral coefficients. Section~\ref{sec:fitting}
details the different model fitting methods,
including Bayesian, maximum likelihood, quasi-maximum likelihood
and exact Whittle likelihood. Our theoretical results are provided
in Section~\ref{sec:theory}. Here, we establish consistency and
asymptotic normality of the proposed estimators.
Section~\ref{sec:sim} illustrates the models effectiveness through a
simulation study, and Section~\ref{sec:conclusion} contains
concluding discussion. Extensions to missing data, imputation, and
signal extraction along with an application of our methodology to
straw yield data from an agricultural experiment, as well as all
proofs, are provided in a Supplementary Appendix (McElroy and Holan \cite{supp}).

\section{The cepstral model and its computation}\label{sec:model}
We begin by introducing some basic concepts about spatial random
fields, and then we specialize to the cepstral random field, with a
focus on computation of autocovariances. References on spatial
random fields include Whittle \cite{whittle1954},
Besag \cite{besag1972correlation}, Rosenblatt \cite{murray1985stationary,rosenblatt2000gaussian}, Solo \cite{solo1986modeling},
Cressie \cite{cressie1993statistics}, Kedem and
Fokianos \cite{kedem2002regression} and
Rue and Held \cite{rueheld:2010}. A random field $\YY= \{ \YY
_{s_1,s_2} \}$ is
a process with indices on a lattice, which in this paper we take to
be $\ZZ^2$. Typically a random field has a mean function
$\mu_{s_1,s_2} = \EE\YY_{s_1,s_2}$, which may be modeled through
regression variables (Cressie \cite{cressie1993statistics}). The mean-corrected field $\YY
- \{ \mu_{s_1,s_2} \}$ will be denoted by $\WW$.

Interest focuses upon weakly stationary random fields, which in
practice is often adequate once mean effects are identified and
accounted for. When all moments are defined, this is equivalent to
the higher cumulants \cite{brillinger2001time} being
dependent only
on lags between the mean-centered variables. The second cumulant
function, or autocovariance function (acf), is defined via
$\operatorname{Cov} ( \YY_{s_1,s_2}, \YY_{r_1,r_2} ) = \EE[ \WW_{s_1,s_2}
\WW_{r_1,r_2} ] = \gamma_{s_1-r_1,s_2-r_2}$ for all $s_1,s_2,r_1,r_2
\in\ZZ$. It is convenient to summarize this second-order structure
through the spectral density $F$ defined on ${[-\pi, \pi]}^2$, which
depends on two frequencies. Letting $Z_j = e^{-i \lambda_j}$ for
$j=1,2$, the spectral density is related to the acf via the formula
%
\begin{equation}
\label{eq:specAcf} F (\lambda_1, \lambda_2 ) = \sum
_{h_1,h_2 \in\ZZ} \gamma_{h_1,h_2} (F) Z_1^{h_1}
Z_2^{h_2}.
\end{equation}
Here we write $\gamma(F)$ for the acf associated with the spectrum
$F$, and it in turn is expressed in terms of $F$ via Fourier
inversion as
%
\begin{equation}
\label{eq:acfCeps} \gamma_{h_1,h_2} (F) = \frac{1}{ 4 \pi^2} \int
_{-\pi}^{\pi} \int_{-\pi
}^{\pi}
F( \lambda_1, \lambda_2 ) Z_1^{-h_1}
Z_2^{-h_2} \,d\lambda_1 \,d\lambda_2.
\end{equation}
As a general notation, let the normalized double integral over both
frequencies be abbreviated by the expression $\langle\cdot
\rangle$, so that $\gamma_{h_1,h_2} (F) = \langle F Z_1^{-h_1}
Z_2^{-h_2} \rangle$ is compactly expressed. Now it follows
elementarily from the commutativity of the field $\YY$ variables
that $\gamma_{h_1,h_2} (F) = \gamma_{-h_1,-h_2} (F)$, and hence the
corresponding $F$ in (\ref{eq:specAcf}) must have mirror
reflectional symmetry through both axes, that is, $F(\lambda_1,
\lambda_2) = F(-\lambda_1, - \lambda_2)$. Furthermore, the acf of a
random field is always positive-definite~\cite{cressie1993statistics} and the corresponding spectrum is
nonnegative \cite{bochner1955harmonic}.

\subsection{The cepstral random field model}

A spatial model for continuous-valued random variables should, at a
minimum, capture second-order structure in the data, which is
summarized through the acf. However, a putative acf may or may not
have nonnegative discrete Fourier transform (DFT)
(\ref{eq:specAcf}), whereas any valid acf of a stationary field
\textit{must} have nonnegative spectrum $F$. One way to ensure our model
has such a valid acf is to model $F$, utilizing some class of
nonnegative functions, and determine the corresponding covariances
via (\ref{eq:acfCeps}). This is the philosophy behind the versatile
exponential time series model of Bloomfield \cite{Bloom:1973}. The
idea there
was to expand the log spectrum in the complex exponential basis
functions, with a truncation of the expansion corresponding to a
postulated model.

The same idea is readily adapted to the spatial context;
Solo \cite{solo1986modeling} seems to be the first formal presentation
of this idea. Given that $F$ is strictly positive and bounded, we
can expand $\log F$ in each frequency concurrently, which yields
\[
\log F ( \lambda_1, \lambda_2) = \sum
_{j_1,j_2 \in\ZZ} \Theta_{j_1,j_2} Z_1^{j_1}
Z_2^{j_2}.
\]
The coefficients $\{ \Theta_{j_1,j_2} = \langle\log F Z_1^{-j_1}
Z_2^{-j_2} \rangle\}$ are called the cepstral coefficients; see
also the recent treatment of Kizilkaya and
Kayran \cite{kizilkaya2005arma}. A pleasing
feature of this representation is that $F^{-1}$ has cepstral
coefficients $\{ - \Theta_{j_1, j_2} \}$. By truncating the
summation, we obtain a parametric model that can approximate the
second-order structure of \textit{any} random field with bounded
spectrum. So we obtain the cepstral model of order $(p_1,p_2)$ given
by
%
\begin{equation}
\label{eq:modelCeps} F ( \lambda_1, \lambda_2) = \exp\Biggl
\{ \sum_{j_1=-p_1}^{p_1} \sum
_{j_2=-p_2}^{p_2} \Theta_{j_1,j_2} Z_1^{j_1}
Z_2^{j_2} \Biggr\}.
\end{equation}
Note that the cepstral coefficient $\Theta_{0,0}$ has no sinusoidal
function multiplying it, and hence $\exp\Theta_{0,0}$ quantifies
the scale of the data. In one dimension, this would be called the
innovation variance; note that $\Theta_{0,0} = \langle\log F
\rangle$. Because the complex exponentials form a complete
orthonormal basis set, it is impossible for two distinct values of
$\Theta$ to produce an identical function $F$; hence the model is
identifiable. Further special cases of the general cepstral field
model are considered in Solo \cite{solo1986modeling}. Because $F$ has
mirror reflectional symmetry, the cepstral coefficients do as well,
that is, $\Theta_{j_1,j_2} = \Theta_{-j_1,-j_2}$.

In order to fit this model to Gaussian data, it is necessary to
compute the acf from a given specification of cepstral coefficients.
We next describe two approaches to this: one is approximate, and the
other is exact. Both differ from the fitting techniques in
Solo \cite{solo1986modeling}, who advocates an asymptotic likelihood (or
Whittle) calculation.

\subsection{Fast calculation of autocovariances}

We here discuss a straightforward discretization of
(\ref{eq:acfCeps}), together with (\ref{eq:modelCeps}),\vadjust{\goodbreak} utilizing
the Riemann approximation. So long as the spectrum is a bounded
function, this method is arbitrarily accurate (since the
practitioner controls the mesh size). In order to accomplish the
computation, without loss of generality let $p_2 = p_1$, so that the
cepstral coefficients are given by a $(2p_1+1) \times(2p_1+1)$ grid
$\Theta$ (if $p_2 < p_1$, just fill in some entries of $\Theta$ with
zeroes).

Now we refer to the entries of $\Theta$ via $\Theta_{j_1,j_2}$ with
$-p_1 \leq j_1,j_2 \leq p_1$, which is a Cartesian mode of indexing;
this differs from the style of indexing pertinent to matrices. We
can map this grid to a matrix $[ \Theta]$ (and back), with the
following rule:
%
\begin{equation}
\label{eq:indexRule} {[ \Theta]}_{k_1,k_2} = \Theta_{k_2-p_1-1,p_1+1-k_1},\qquad
\Theta_{j_1, j_2} = {[ \Theta]}_{p_1+1-j_2, j_1+p_1+1}
\end{equation}
for $1 \leq k_1,k_2 \leq2p_1+1$ and $-p_1 \leq j_1,j_2 \leq p_1$.
We will consider a set of frequencies $\{ \ell_1 \pi/ M, \ell_2 \pi
/ M
\}$ for $-M
\leq\ell_1, \ell_2 \leq M$, which is an order $M$ discretization of
${[-\pi,\pi]}^2$. Suppose that we wish to compute
the grid of autocovariances given by $\Gamma= { \{
\gamma_{h_1,h_2} \} }_{h_1,h_2=-H}^H$ for some maximal lag $H$. To
that end, we consider a complex-valued $2p_1+1 \times2M+1$ matrix $E$ with
entries $E_{k_1, k_2} = \exp\{ i \pi(p_1+1-k_1)(M-k_2+1) M^{-1} \}$ for
$k_1=1,2,\ldots, 2p_1+1$ and $k_2 = 1,2, \ldots, 2M+1$, and also
define a $(2H+1) \times(2M+1)$ dimensional matrix $G$ via
$G_{j_1, j_2} = \exp\{ i \pi(H+1-j_1)(M+1-j_2) M^{-1} \} $.
Then with $[\Gamma]$ defined via ${[ \Gamma] }_{k_1, k_2} = \gamma
_{k_2-H-1, H+1 -
k_1}$, the formula
%
\begin{equation}
\label{eq:acfApprox} [ \Gamma] \approxeq{(2M+1)}^{-2} G \exp\bigl\{
\overline{E}^{\prime} [\Theta] E \bigr\} \overline{G}^{\prime}
\end{equation}
provides a practical method of computation. In this formula, which
is derived in the Supplement's Appendix B, we have written the
exponential of a matrix, which
here is \textit{not} the ``matrix exponential,'' but rather just
consists of exponentiating each entry of the matrix.
So (\ref{eq:acfApprox}) produces an arbitrarily fine approximation
to the acf (taking $M$ as large as desired). The algorithm takes a
given $\Theta$, produces $[\Theta]$ via (\ref{eq:indexRule}),
computes $E$ and $G$ (ahead of time, as they do not depend upon the
parameters) and determines $[\Gamma]$ via~(\ref{eq:acfApprox}).

\subsection{Exact calculation of autocovariances}

We now present an exact\break  method for computing the acf from the
cepstral matrix. Our approach is similar to that of Section~3 of
Kizilkaya and
Kayran \cite{kizilkaya2005arma}, though with one important difference.
They present an algorithm for computing cepstral coefficients from
known coefficients of an ARMA random field. Instead, we take the
cepstral coefficients as given, compute coefficients of certain
corresponding MA random fields, and from there obtain the acf. In
order to fit the Gaussian likelihood, we need to compute the acf
from the cepstral matrix, not the reverse.

We introduce the device of a ``causal'' field and a ``skew'' field as
follows. The causal field is an MA field that only involves
coefficients with indices in the positive quadrant, whereas the
skew\vadjust{\goodbreak}
field essentially is defined over the second quadrant. More
precisely, we have
%
\begin{eqnarray}
\label{eq:causalField} \gamma_{s_1,s_2} (\Psi) & = &\sum
_{k_1,k_2 \geq0 } \psi_{s_1+k_1,s_2+k_2} \psi_{k_1,k_2},
\nonumber
\\[-8pt]
\\[-8pt]
\nonumber
{ \biggl\llvert\sum_{j_1,j_2 \geq0 } \psi_{j_1,j_2}
Z_1^{j_1} Z_2^{j_2} \biggr\rrvert
}^2 & = &\sum_{s_1,s_2 \in\ZZ} \gamma_{s_1,s_2}
(\Psi) Z_1^{s_1} Z_2^{s_2}
\end{eqnarray}
for the causal field. The causal field may be written formally (in
terms of backshift operators $B_1, B_2$) as $\Psi(B_1, B_2) =
\sum_{j_1,j_2 \geq0 } \psi_{j_1,j_2} B_1^{j_1} B_2^{j_2}$.
That is, the $\psi_{j_1,j_2}$ coefficients define the moving average
representation of the causal field, and $\{ \gamma_{s_1,s_2} (\Psi)
\}$ is its acf. It is important that we set $\psi_{0,0} = 1$.
Similarly, let $\Phi(B_1, B_2) = \sum_{j_1,j_2 \geq0 }
\phi_{j_1,j_2} B_1^{-j_1} B_2^{j_2}$ for the skew-field, which
in the first index depends on the forward shift operator $B_1^{-1}$,
but on the backshift operator $B_2$ in the second index. Thus
%
\begin{eqnarray}
\label{eq:skewField} \gamma_{s_1,s_2} (\Phi) & =& \sum
_{k_1,k_2 \geq0 } \phi_{s_1+k_1,s_2+k_2} \phi_{k_1,k_2},
\nonumber
\\[-8pt]
\\[-8pt]
\nonumber
{ \biggl\llvert\sum_{j_1,j_2 \geq0 } \phi_{j_1,j_2}
Z_1^{-j_1} Z_2^{j_2} \biggr\rrvert
}^2 & =& \sum_{s_1,s_2 \in\ZZ} \gamma_{s_1,s_2}
(\Phi) Z_1^{-s_1} Z_2^{s_2}.
\end{eqnarray}
We also have two time series, corresponding to the axes of the
cepstral matrix, given by $\Xi(B_1) = \sum_{j_1 \geq0 } \xi_{j_1}
B_1^{j_1}$ and $\Omega(B_2) = \sum_{j_2 \geq0} \omega_{j_2}
B_2^{j_2}$, which have acfs $\gamma_{h_1} (\Xi) = \sum_{k_1 \geq0}
\xi_{k_1+h_1} \xi_{k_1} $ and $\gamma_{h_2} (\Omega) = \sum_{k_2
\geq0} \omega_{k_2+h_2}, \omega_{k_2} $, respectively. Now each of
these MA random fields has a natural cepstral representation, such
that their acfs can be combined to produce the cepstral acf, as
shown in the following result.

%
\begin{Proposition}
\label{prop:cepstral2MA}
The acf of the cepstral model is given by
%
\begin{eqnarray}
\label{eq:acfExpress}\qquad&& \gamma_{h_1,h_2} (F)
\nonumber
\\[-4pt]
\\[-12pt]
\nonumber
&&\qquad= e^{\Theta_{0,0}} \sum
_{j_1,j_2 \in\ZZ} \gamma_{j_1,j_2} (\Phi) \biggl[ \sum
_{k_1,k_2
\in\ZZ} \gamma_{h_1+j_1-k_1,h_2-j_2-k_2} (\Psi) \gamma_{k_1} (
\Xi) \gamma_{k_2} (\Omega) \biggr],
\end{eqnarray}
where $\gamma(\Phi)$, $\gamma(\Psi)$, $\gamma(\Xi)$ and $\gamma
(\Omega)$ can be calculated in terms of their coefficients, which
are recursively given by
%
\begin{eqnarray}
\label{eq:betaRecurs} \psi_{j_1,j_2} & =& \frac{1}{j_1} \sum
_{k_1 = 1 }^{p_1} k_1 \Biggl( \sum
_{k_2=1}^{j_2} \psi_{j_1-k_1,j_2-k_2}
\Theta_{k_1,k_2} \Biggr),
\\
\label{eq:rhoRecurs} \phi_{j_1,j_2} & =& \frac{1}{j_1} \sum
_{k_1 =1 }^{p_1} k_1 \Biggl( \sum
_{k_2=1}^{j_2} \phi_{j_1-k_1,j_2-k_2}
\Theta_{-k_1,k_2} \Biggr),
\\
\label{eq:deltaRecurs} \xi_{j_1} & =& \frac{1}{j_1} \sum
_{k_1 =1 }^{p_1} k_1 \Theta_{k_1,0}
\xi_{j_1-k_1},
\\
\label{eq:etaRecurs} \omega_{j_2} & = &\frac{1}{j_2} \sum
_{k_2 =1}^{p_1} k_2 \Theta_{0,k_2}
\omega_{j_2-k_2}
\end{eqnarray}
for $j_1 \geq1$ and $j_2 \geq1$.
\end{Proposition}

Proposition \ref{prop:cepstral2MA} gives recursive formulas. In the
causal case, one would compute $ \psi_{1,1}, \psi_{2,1}, \ldots,
\psi_{p_1,1}, \psi_{1,2}, \psi_{2,2}, \ldots,$ etc. Alternative
computational patterns could be utilized, noting that
$\psi_{j_1,j_2}$ only requires knowledge of $\psi_{\ell_1, \ell_2}$
with $\ell_1 < j_1$ and $\ell_2 < j_2$.
When $p_1 = \infty$, equation (\ref{eq:acfExpress}) gives the
precise mapping of cepstral coefficients to various MA coefficients,
and ultimately to the autocovariance function. If $p_1 < \infty$,
it provides an algorithm for determining autocovariances for a given
cepstral model. These formulas are already much more complicated
than in the time series case (see \cite{pourahmadi1984}), and for
higher dimensional fields become intractable.

\section{Model fitting methods}\label{sec:fitting}
In this section we give additional details on various methods for
fitting cepstral random field models and present some tools for
refining specified models. Once a model is specified, we can
estimate the parameters via exact maximum likelihood, Bayesian
posterior simulation, an approximate Whittle likelihood or an exact
Whittle likelihood. We focus on these four techniques due to their
mixture of being flexible and possessing good statistical
properties.

We first define Kullback--Leibler (KL) discrepancy, the exact Whittle
likelihood and the quasi-maximum likelihood estimate (QMLE), and
then we proceed to describe the distributional behavior of the
maximum likelihood estimates (MLEs) and QMLEs, extending the results
of Mardia and Marshall \cite{mardia1984maximum} to non-Gaussian
fields, under an
expanding domain asymptotic theory. These results, proved for
fairly general linear random fields with regression effects, are
then specialized to the case of the cepstral field, and model
selection is afterwards described.

\subsection{Random field data}
Now we proceed to discuss spatial modeling (here we do not assume a
cepstral random field structure), adapting the vector time series
treatment in Taniguchi and
Kakizawa \cite{taniguchi2000asymptotic}. Suppose that our data
comes to us in gridded form, corresponding to a $N_1 \times N_2$
matrix $\YY^N$ (with $\WW^N$ denoting the de-meaned version). We
use the notation $N = \sqrt{N_1 \cdot N_2}$, so that $N^2$ is the
sample size. Both $\YY^N$ and $\WW^N$ can be vectorized into length
$N^2$ vectors $Y$ and $W$ via the so-called lexicographical rule
\[
Y_k = \YY^N_{s_1,s_2},\qquad W_k =
\WW^N_{s_1,s_2},\qquad k = N_2 (s_1-1) +
s_2.
\]
Here, $Y = \operatorname{vec} ( { \YY^N }^{\prime} )$, where $\operatorname{vec}$
stands for the vector operation on a matrix, and $\prime$ is the
transpose. Note that $s_1 -1 = k \operatorname{div} N_2$ and $s_2 = k
\operatorname{mod} N_2$. Also let $\mu= \EE Y$, so that $\mu_k = \EE
Y_k = \EE\YY^N_{s_1,s_2} = \mu_{s_1,s_2}$. In the simplest scenario
the mean matrix $\{ \mu_{s_1,s_2} \}$ is constant with respect to
the indices $s_1,s_2$. More generally, we might model the mean
through regressor functions defined upon the grid, that is,
$\mu_{s_1,s_2} = \sum_{\ell= 1}^L \beta_{\ell} X_{\ell} (s_1,s_2) $
for some specified lattice functions ${ \{ X_{\ell} \}
}_{\ell=1}^L$. Then
\[
\mu_k = \sum_{\ell=1}^L
\beta_{\ell} X_{\ell} ( k \operatorname{div} N_2 + 1, k
\operatorname{mod} N_2) = \sum_{\ell=1}^L
\beta_{\ell} \widetilde{X}_{\ell} ( k )
\]
maps each $X_{\ell}$ from a lattice function to a function $\widetilde
{X}_{\ell}$ of the
natural numbers. The parameters $\beta_1, \beta_2, \ldots,
\beta_L$ then enter the regression linearly, and we can express
things compactly via $\mu= \widetilde{X} \beta$, where
$\widetilde{X}$ is the regression matrix with columns given by the
various $\widetilde{X}_{\ell}$.

The spectral density of a mean zero random field $\WW$ has already
been defined in~(\ref{eq:specAcf}), and the DFT of the field is now
defined as
\[
\widetilde{\WW} (\lambda_1, \lambda_2) = \sum
_{t_1=1}^{N_1} \sum_{t_2 =
1}^{N_2}
\WW_{t_1, t_2}^N e^{-i \lambda_1 t_1 } e^{-i \lambda_2 t_2 } = \sum
_{t_1 = 1}^{N_1} \sum
_{ t_2 = 1}^{N_2} W_{N_2 (t_1 -1) + t_2} Z_1^{t_1}
Z_2^{t_2}
\]
for $\lambda_1, \lambda_2 \in[- \pi, \pi]$. Note that we define
this DFT over all pairs of frequencies, not just at the so-called
Fourier frequencies. Also the DFT depends on $\beta$ through the
mean-centering; if we center the data $\YY^N$ by any regression
parameter other than the true $\beta$, denoted $\widetilde{\beta}$,
some bias will be introduced. The periodogram will be defined at all
frequencies and is proportional to the squared magnitude of the
DFT,
\[
I_{\widetilde{\beta}} (\lambda_1, \lambda_2) =
N^{-2}  \bigl| \widetilde{\WW} (\lambda_1,
\lambda_2)\bigr |^2 = \sum_{|h_1|<N_1}
\sum_{ |h_2 | < N_2} \gamma_{h_1,h_2} (I_{\widetilde{\beta}})
Z_1^{h_1} Z_2^{h_2}.
\]
Here $\gamma_{h_1,h_2} (I_{\beta}) $ is defined as the
sample acf of the series demeaned by $\mu= \widetilde{X} \beta$ (see
Supplementary Appendix B for more detail);
moreover it satisfies (\ref{eq:acfCeps}) with $F$ replaced by
$I_{\beta}$. We also will consider an
unbiased acf estimate given by
\[
\widehat{\gamma}_{h_1, h_2} (I_{\beta})= \frac{ N^2 }{ (N_1-|h_1|)(N_2
-|h_2|) }
\gamma_{h_1, h_2} (I_{\beta}).
\]
We emphasize that the computation of this periodogram requires a
choice of $\beta$, and so is written $I_{\beta}$. This can be used
to assess the frequency domain information in the random field along
any row or column; the periodogram can also be viewed as a crude
estimate of the spectral density $F$ \cite{cressie1993statistics}.

In our context the treatment of the periodogram differs from the
treatment provided in Fuentes \cite{fuentes2002spectral}. In particular,
we consider the periodogram defined at all frequencies, not just
the\vadjust{\goodbreak}
Fourier frequencies. Additionally, the asymptotic properties
developed in Fuentes \cite{fuentes2002spectral} rely on shrinking domain
asymptotics, whereas our asymptotic arguments rely on an expanding
domain. Finally, our periodogram is defined in terms of a mean
centered random field and, thus, explicitly depends on the
regression parameters $\beta$.

\subsection{Model fitting criteria}
Let the covariance matrix of $\WW^N$ be denoted
$\Sigma(\widetilde{F})$, which is defined via $\Sigma
(\widetilde{F}) = \EE W W^{\prime}$; the resulting block-Toeplitz
structure of this matrix is analyzed in Section~\ref{sec:theory}.
The entries of this matrix can be determined from $\widetilde{F}$
via the algorithms of Section~\ref{sec:model}, along with careful
bookkeeping. A model for the data involves a spectrum
$F_{\theta}$---let the associated block-Toeplitz covariance matrix be denoted
$\Sigma(F_{\theta})$---which is hoped to be a suitable
approximation to $\Sigma(\widetilde{F})$. Model fitting can be
performed and assessed through the Kullback--Leibler (KL)
discrepancy, just as with time series. Although KL is mentioned in
Solo \cite{solo1986modeling} and Cressie \cite
{cressie1993statistics}, we
provide an in-depth treatment here; see Lemma \ref{lem:Fourier},
for example. If $F$ and $G$ are two (mean zero) random field
spectral densities, their KL discrepancy is defined to be
\[
\operatorname{KL} (F, G) = \langle\log F + G/ F \rangle.
\]
This is a convenient mechanism, since KL is convex in $F$. As
$\beta$ parametrizes mean effects, we let $\theta$ be a parameter
vector describing the second-order structure. If the true data
process\vspace*{1pt} has spectrum $\widetilde{F}$, and we utilize a model with
spectrum $F_{\theta}$, then KL($F_{\theta}$, $\widetilde{F}$) can be
used to assess proximity of the model to truth. The convexity of KL
guarantees that when the model is correctly specified, the true
parameter $\widetilde{\theta}$ minimizes the discrepancy. When the
model is misspecified, the minima $\widetilde{\theta}$ are called
pseudo-true values (cf.~\cite{taniguchi2000asymptotic}). For the
cepstral model, the parameter vector is $\theta= J \operatorname{vec}
\Theta$, where $J$ is a selection matrix that eliminates
redundancies in $\Theta$ due to symmetry. The full parameter vector
is written $\phi$, where $\phi^{\prime} = {[ \theta^{\prime},
\beta^{\prime} ]}$.

It is natural to use KL to fit models as well. For this, consider
KL($F_{\theta}$, $I_{\beta}$)---which is called the exact Whittle
likelihood---and minimize with respect to $\theta$, which produces
by definition the estimate $\widehat{\theta}_{\mathrm{QMLE}}$. Then using
(\ref{eq:acfCeps}) we obtain the practical expression
\[
\operatorname{KL} (F_{\theta}, I_{\beta} ) = \langle\log
F_{\theta} \rangle+ \sum_{|h_1 | < N_1} \sum
_{ |h_2| < N_2} \gamma_{h_1,
h_2} (I_{\beta}) \cdot
\gamma_{h_1, h_2} \bigl(F_{\theta}^{-1}\bigr).
\]
This assumes that the correct regression parameters have been
specified. In the case that $F_{\theta}$ is a cepstral spectrum
(\ref{eq:modelCeps}), the above expression is even easier to
compute: $\langle\log F_{\theta} \rangle = \Theta_{0,0}$ and
$\gamma(F_{\theta}^{-1}) = \gamma(F_{- \theta})$, that is, multiply
each cepstral coefficient by $-1$ to obtain the acf of
$F_{\theta}^{-1}$ from the acf of $F_{\theta}$.

Unfortunately, $\gamma_{h_1, h_2} (I_{\widetilde{\beta}})$ is biased
as an estimate of
$\gamma_{h_1, h_2} (\widetilde{F})$, and this has a nontrivial
impact for spatial data, though not for time series. Essentially,
the presence of ``corners'' in the observed data set reduces the
number of data points that are separated by a given lag $(h_1,
h_2)$; if either of $|h_1|$ or $|h_2|$ is large, we have a very
biased estimate. Note, the impact of corners can be visualized by
comparing the volume
of a $d$-dimensional cube with that of an inscribed ball; the ratio is
$\pi/(2d)$ for $d \geq2$,
which tends to zero as $d$ increases. Thus, corners increasingly
dominate the region as
$d$ increases, which interferes with one's ability to measure
correlation as
a function of lag. This effect is more pronounced as the dimension
increases. For this reason,\vspace*{1.5pt} we propose using
$\widehat{\gamma}_{h_1, h_2} (I_{\beta})$ instead of $\gamma_{h_1,
h_2} (I_{\beta})$,
because $\EE\widehat{\gamma}_{h_1, h_2} (I_{\widetilde{\beta}}) =
\gamma_{h_1, h_2}
(\widetilde{F})$. Let us call the modified KL($F_{\theta}$,
$I_{\beta}$) by $\widehat{\operatorname{KL}} (F_{\theta})$
\[
\widehat{\operatorname{KL}} (F_{\theta} ) = \langle\log F_{\theta}
\rangle+ \sum_{|h_1| < N_1} \sum
_{|h_2| < N_2 } \widehat{\gamma}_{h_1, h_2} (I_{\beta})
\cdot\gamma_{h_1, h_2} \bigl(F_{\theta}^{-1} \bigr).
\]
Using this criterion instead will produce asymptotically normal
cepstral parameter estimates, and therefore is to be preferred.

A drawback of utilizing $\widehat{\gamma}_{h_1,h_2} (I_{\beta})$ is
that the
corresponding spectral esti\-mate---the DFT of the
unbiased sample acf---need not be positive at all frequencies.
Although this is irrelevant asymptotically, in finite samples it can
interfere with inference. The time domain representation of
$\widehat{\operatorname{KL}} (F_{\theta})$ can still be computed, of course but
the second term in its formula might not be positive. Other types
of autocovariance estimators could be utilized, being based on
other kinds of spectral estimators (see Politis and Romano \cite
{politis1996flat}
for a discussion of the tradeoff between bias and nonnegativity of
the spectral estimate). These alternative estimators might be
based on convolving the periodogram with a spectral window, or
equivalently by using a taper (or lag window) with the sample acf.
Tapers are known to modify the bias and variance properties of
spectral estimators in time series; see Guyon \cite{guyon1982parameter},
Dahlhaus and
K{\"u}nsch \cite{dahlhaus1987edge} and Politis and Romano \cite
{politis1995bias}.

If even faster computation of the objective function is desired, we
may discretize $\widehat{\operatorname{KL}}$ and utilize values of $F$
directly, without having to compute the inverse DFT $\gamma
(F_{\theta}^{-1})$. The result is the approximate Whittle
likelihood, denoted $\widehat{\operatorname{KL}}_N$, and is obtained by
discretizing the integrals in KL($F_{\theta}$, $I_{\beta}$) with a
mesh corresponding to Fourier frequencies, but replacing $I_{\beta}$
with the DFT of the $\widehat{\gamma}_{h_1, h_2} (I_{\beta})$
sequence, denoted by $\widehat{I}_{\beta}$. Then the discrepancy is
\[
\label{eq:KLunbiased} \widehat{\operatorname{KL}}_N (F_{\theta}) =
N^{-2} \sum_{j_1 = -N_1}^{N_1} \sum
_{j_2 = -N_2}^{N_2} \biggl\{ \log F_{\theta}
\biggl( \frac{\pi j_1}{ N_1}, \frac{\pi
j_2}{N_2} \biggr) + \frac{ \widehat{I}_{\beta} ( { \pi j_1
}/{N_1}, { \pi
j_2 }/{ N_2} ) }{
F_{\theta} ( { \pi j_1}/{N_1}, { \pi j_2 }/{ N_2} )
} \biggr\},
\]
which can be minimized with respect to $\theta$. The resulting estimate
has asymptotic properties identical to the QMLE, and in practice one
may use either $\widehat{\operatorname{KL}}$ or $\widehat{\operatorname{KL}}_N$
according to computational convenience. It will be convenient to
present a notation for this double discrete sum, which is a Fourier
approximation to $\langle\cdot\rangle$, denoted by ${\langle\cdot
\rangle}_N$; then $\widehat{\operatorname{KL}}_N (F_{\theta}) = {\langle
\log
F_{\theta} +
\widehat{I}_{\beta}/ F_{\theta} \rangle}_N$.

We can also extend the KL formula to handle regression effects,
%
\begin{equation}
\label{eq:KLquadratic} \operatorname{KL} (F_{\theta}, I_{\beta} ) = \langle
\log
F_{\theta} \rangle+ N^{-2} {(Y - \widetilde{X}
\beta)}^{\prime} \Sigma\bigl(F^{-1}_{\theta
}\bigr) {(Y -
\widetilde{X} \beta)}.
\end{equation}
This formula is proved in Supplementary Appendix B. We propose using (\ref
{eq:KLquadratic}) to
estimate regression parameters, but $\theta$ is to be determined by
$\widehat{\operatorname{KL}}$. The formula for the regression QMLE is then
%
\begin{equation}
\label{eq:betaQMLE} \widehat{\beta}_{\mathrm{QMLE}} = { \bigl[ \widetilde
{X}^{\prime}
\Sigma\bigl(F_{\widehat{\theta}_{\mathrm{QMLE}}}^{-1} \bigr) \widetilde{X}
\bigr]}^{-1} \widetilde{X}^{\prime} \Sigma\bigl(F_{\widehat
{\theta}_{\mathrm{QMLE}}}^{-1}
\bigr) Y,
\end{equation}
where $\widehat{\theta}_{\mathrm{QMLE}}$ minimizes $\widehat{\operatorname{KL}}
(F_{\theta})$, which in turn depends upon $\widehat{\beta}_{\mathrm{QMLE}}$
through $\widehat{I}_{\beta}$. These formulas do not apply when we use
the approximate Whittle,
although the same asymptotic properties will hold as for the exact
Whittle.

On the other hand, we can also compute the exact Gaussian likelihood
for the field. The log Gaussian likelihood is equal (up to
constants) to
%
\begin{equation}
\label{eq:GaussLik} \mathcal{L} (\theta, \beta) = -\tfrac{1}{2}
\log\bigl|
\Sigma(F_{\theta
})\bigr | - \tfrac{1}{2} {(Y - \widetilde{X}
\beta)}^{\prime} \Sigma^{-1} (F_{\theta}) {(Y -
\widetilde{X} \beta)}.
\end{equation}
Maximizing this function with respect to $\theta$ yields the MLE
$\widehat{\theta}_{\mathrm{MLE}}$; also $\widehat{\beta}_{\mathrm{MLE}}$ is given by
the generalized least squares (GLS) estimate by standard arguments
(see~\cite{mardia1984maximum}),
%
\begin{equation}
\label{eq:betaMLE} \widehat{\beta}_{\mathrm{MLE}} = { \bigl[ \widetilde
{X}^{\prime}
\Sigma^{-1} (F_{\widehat{\theta}_{\mathrm{MLE}}}) \widetilde{X} \bigr]}^{-1}
\widetilde{X}^{\prime} \Sigma^{-1} (F_{\widehat{\theta}_{\mathrm{MLE}}}) Y,
\end{equation}
which expresses the regression parameter in terms of
$\widehat{\theta}_{\mathrm{MLE}}$. For the computation of
(\ref{eq:GaussLik}) we must calculate the acf corresponding to
$F_{\theta}$, which can be done using the algorithms of
Section~\ref{sec:model}. Contrast (\ref{eq:GaussLik}) with
(\ref{eq:KLquadratic}); they are similar, the main difference being
the replacement of the inverse of $\Sigma(F_{\theta})$ by $\Sigma
(F^{-1}_{\theta})$, which is equal to $\Sigma( F_{- \theta} )$ for
the cepstral model.

Most prior literature on random fields seems to utilize
approximate Whittle estimation, or QMLE, since the objective
function is quite simple to write down.
The parameter MLEs do not have the bias problem of QMLEs, discussed
above, but require more effort to compute due to matrix inversion.
We can use the approximate algorithm given by equation
(\ref{eq:acfApprox}), together with (\ref{eq:GaussLik}), to compute
the MLEs. The QMLEs, based on unbiased $\widehat{\gamma}_{h_1,h_2}
(I_{\beta})$ acf estimates, are faster to compute than MLEs and
enjoy the same asymptotic normality and efficiency.

However, if one prefers a Bayesian estimation of $\theta$ (and
$\beta$), it is necessary to compute $\exp\mathcal{L} (\theta,
\beta)$, which is proportional to the data likelihood $p (Y \vert
\theta, \beta)$. The posterior for $\theta$ is proportional to the
likelihood times the prior, and one can use Markov chain Monte Carlo
(MCMC) methods to approximate $p(\theta\vert Y)$
\cite{geweke2005contemporary}. Recall that the mean of this
distribution, which is the conditional expectation of $\theta$ given~$Y$, is called the posterior mean, and will be denoted
$\widehat{\theta}_B$.

Note that equations (\ref{eq:GaussLik}) and (\ref{eq:betaMLE}) can
each be used in an iterative estimation scheme. To determine the
MLE, minimize (\ref{eq:GaussLik}) to obtain an estimate of $\theta$
for a given $\beta$ computed via (\ref{eq:betaMLE}); then update
$\widehat{\beta}$ by plugging into (\ref{eq:betaMLE}), and iterate.
For the QMLE, de-mean the data by computing $Y- X \beta$ and
determining $\widehat{I}_{\beta}$, for a given $\beta$, and then
minimize either $\widehat{\operatorname{KL}}$ or $\widehat{\operatorname{KL}}_N$ to
obtain $\theta$ estimates (either exact or approximate); then update
$\beta$ by plugging into (\ref{eq:betaQMLE}) and iterate. From now
on, we refer to these estimates as the exact/approximate QMLEs [if
using biased acf estimates $\gamma_{h_1,h_2} (I_{\beta})$, only
consistency holds, and not asymptotic normality].

\subsection{Distributional properties of parameter estimates}

We now provide a description of the asymptotics for the various
estimates; a rigorous treatment is given in
Section~\ref{sec:theory}, with formal statements of sufficient
conditions and auxiliary results. First, the Bayesian estimates
$\widehat{\theta}_B$ and $\widehat{\beta}_B$ are consistent when the
data is a Gaussian random field that satisfies suitable regularity
conditions (Theorem \ref{thmm:Bayes}). For the frequentist case,
recall that the number of observations equals $N_1 \cdot N_2$, so
that a central limit theorem result requires scaling by $\sqrt{N_1
\cdot N_2}$; we require that both dimensions expand, that is, $\min\{
N_1, N_2 \} \tends\infty$. Let the Hessian of the KL be denoted
$H(\theta) = \nabla\nabla^{\prime} \operatorname{KL} (F_{\theta}, \widetilde{F})$,
which will be invertible at the unique pseudo-true value
$\widetilde{\theta}$ by assumption. Then the exact QMLE, approximate
QMLE and MLE for $\theta$ are all consistent, and are also
asymptotically normal at rate $N$ with mean $\widetilde{\theta}$ and
variance $H^{-1} (\widetilde{\theta}) V (\widetilde{\theta})
H^{\dag} (\widetilde{\theta})$, where $\dag$ denotes inverse
transpose and $V(\theta) = 2 \langle\widetilde{F}^2 \nabla
F_{\theta}^{-1} \nabla^{\prime} F_{\theta}^{-1} \rangle$. (This
assumes that the fourth cumulants are zero; otherwise a more
complicated expression for $V$ results, involving the fourth-order
spectral density.) The estimates of the regression parameters are
asymptotically normal and independent of the $\theta$ estimates
(when the third cumulants are zero), for all three types of
estimates.

These theoretical results can be used to refine models. Typically,
one uses these types of asymptotic results under the null hypothesis
that the model is correctly specified, so that $\widetilde{\theta}$
is the true parameter and $V = 2 \langle\nabla\log F_{\theta}
\nabla^{\prime} \log F_{\theta} \rangle$, which equals twice $H$.
See McElroy and Holan \cite{mcelroy2009local} and McElroy and
Findley \cite{mcelroy2010selection} for
more exposition on model misspecification in the frequency domain.
Thus, the asymptotic variance is twice the inverse Hessian, or the
inverse of $H/2$. Note that the Fisher information matrix is the
Hessian of the asymptotic form of the Whittle likelihood, and hence
is equal to one half of the Hessian of KL, that is, $H/2$. Therefore
when the model is correctly specified, parameter estimation is
efficient.

Furthermore, the Fisher information matrix has a particularly
elegant form in the case of a cepstral model. The gradient of the
log spectrum is in this case just the various $Z_1^{j_1}$ or
$Z_2^{j_2}$, so that as in the time series case the Hessian equals
twice the identity matrix (because of mirror reflectional symmetry
in $\Theta$, there is a doubling that occurs), except for the case
of the entry corresponding to $\Theta_{0,0}$---in this case the
derivative of the log spectrum with respect to $\Theta_{0,0}$ equals
one. Thus the Fisher information matrix for all the parameters
except $\Theta_{0,0}$ is equal to the identity matrix, and hence the
asymptotic variance of any cepstral coefficient estimate is $N^{-2}$
(or $2 N^{-2}$ in the case of $\Theta_{0,0}$). The lack of
cross-correlation in the parameter estimates asymptotically
indicates there is no redundancy in the information they convey,
which is a type of ``maximal efficiency'' in the cepstral model.

In terms of model-building with cepstral random fields, one
procedure is the following: postulate a low order cepstral field
model (e.g., order $p_1=1$) and jointly test for whether any
coefficients (estimated via MLE or QMLE) are equal to zero. We might consider
expanding the model---in the direction of one spatial axis or
another as appropriate---if coefficients are significantly
different from zero. Although this is not an \textit{optimal} method of
model selection, this type of forward addition strategy would
stop once all additional coefficients are negligible. Alternatively,
one could start with a somewhat larger cepstral model, and
iteratively delete insignificant coefficients.

Gaussian likelihood ratio test statistics can be utilized for nested
cepstral models, along the lines given in
Taniguchi and
Kakizawa \cite{taniguchi2000asymptotic}---which ultimately just depend on
the asymptotic normality of the parameter estimates---in order to
handle batches of parameters concurrently. Model selection and
assessment can also be assisted by examination of spatial residuals,
which are defined by applying the inverse square root of the
estimated data covariance matrix $\Sigma(F_{\widehat{\theta}})$ to
the vectorized centered data $W$---the result is a vectorized
residual sequence, which should behave like white noise if the model
has extracted all correlation structure. Note that examining
whiteness of the vectorized residuals is equivalent to looking at
all spatial correlations of the spatial residuals defined by undoing
the vec operation. In the context of lattice data, one popular
method for testing the null hypothesis of the absence of spatial
autocorrelation is through the use of Moran's $I$ statistic
\cite{moran1950notes,cliff1981spatial}.
For a
comprehensive discussion regarding Moran's $I$ statistic and its
limitations see Cressie \cite{cressie1993statistics},
Li, Calder and Cressie \cite{li2007beyond},
Cressie and Wikle \cite{cressiewikle2011} and the
references therein. In our case (Supplementary Appendix A), we will evaluate
goodness-of-fit by applying Moran's $I$
statistic to the spatial residuals obtained from the estimated
model.

\section{Theory of inference}\label{sec:theory}

This section provides rigorous mathematical results regarding the
inference problems delineated in Section~\ref{sec:fitting}. We do
not assume a cepstral random field process, retaining greater
generality, but assume a fair amount of regularity on the higher
moments of the field through the Brillinger-type cumulant conditions
\cite{brillinger2001time}. We need not assume the field is
Gaussian for Theorem \ref{thmm:AsympFreq}, but we require a Gaussian
assumption for Theorem \ref{thmm:Bayes}. We first list technical
assumptions, and then describe the mathematical results.

Previous rigorous work on asymptotics for parameter estimates of
lattice random fields includes Guyon \cite{guyon1982parameter} and
Mardia and Marshall \cite{mardia1984maximum}. Although Solo~\cite
{solo1986modeling}
advocates the approximate QMLE method in practice, asymptotic
results are not proved in that paper. Our approach, like
Mardia and Marshall \cite{mardia1984maximum} and Pierce~\cite
{pierce1971}, handles regression
effects together with parameter estimates, but we utilize broader
data process assumptions formulated in terms of cumulants;
Mardia and Marshall~\cite{mardia1984maximum} assumes that the random
field is Gaussian,
whereas we do not. Pierce~\cite{pierce1971} treats the $d=1$
time series
case, allows for non-Gaussian marginals and shows that skewness
can produce asymptotic correlation between regression and model
parameter estimates; an analogous story for $d=2$ is described in
Theorem \ref{thmm:AsympFreq}. Our contribution broadens the
applicability of Mardia and Marshall \cite{mardia1984maximum} to
non-Gaussian fields,
and we moreover provide sufficient conditions under which our Lemma
\ref{lem:toeplitz} yields the validity of condition (iii) of Theorem~2 of Mardia and Marshall \cite{mardia1984maximum}. This highlights our
frequentist
contribution; for Bayesian analysis, we are unaware of any published
work on asymptotic concentration for random fields. Theorem
\ref{thmm:Bayes} assumes a Gaussian field, which is natural given
that the likelihood is Gaussian.

\subsection{Regularity assumptions}

We first set out some notation and working assumptions: define a
block-Toeplitz matrix $\Sigma(F)$ associated with spectral density
$F$ and an $N_1 \times N_2$ data matrix $\YY^N$ to be $N^2 \times
N^2$, with $j_1,k_1$th block (for $1 \leq j_1,k_1 \leq N_1$) given
by the $N_2 \times N_2$-dimensional matrix $\Sigma(F_{j_1-k_1})$,
which is defined as follows. If we integrate over the second
variable of $F$ we obtain a function of the first frequency,
\[
F_{h_1} ( \lambda_2) = \frac{1}{ 2 \pi} \int
_{-\pi}^{\pi} F(\lambda_1,
\lambda_2) e^{i h_1 \lambda_1} \,d\lambda_1\qquad \mbox{with } 0
\leq h_1 < N_1.
\]
Then $\Sigma(F_{h_1})$ is the $N_2 \times N_2$-dimensional matrix of
inverse DFTs of $F_{h_1}$,
with $j_2,k_2$th entry given by
\begin{eqnarray*}
\gamma_{h_1,j_2-k_2} (F)& =& \frac{1}{ 4 \pi^2 } \int_{-\pi}^{\pi}
\int_{-\pi}^{\pi} F(\lambda_1,
\lambda_2) e^{i h_1 \lambda_1} \,d\lambda_1 e^{i (j_2-k_2)
\lambda_2} \,d
\lambda_2 \\
&=& \bigl\langle F Z_1^{-h_1}
Z_2^{k_2-j_2} \bigr\rangle.
\end{eqnarray*}
Based on how we have defined $\WW^N$ and $W = \operatorname{vec} (
{\WW^N}^{\prime})$, it follows that $\Sigma(\widetilde{F}) = \EE[W
W^{\prime}]$, where $\widetilde{F}$ corresponds to the true data
process. That is, lexicographical ordering of a stationary field
produces this structure in the covariance matrix; there are $N_1^2$
blocks, each of which are $N_2 \times N_2$-dimensional.

Also let $\mathcal{F}$ denote the set of admissible spectra for
two-dimensional random fields, defined as follows. For any spatial
autocovariance function $\{ \gamma_{h_1,h_2} \}$, consider the sums
$S_{h_1, \cdot} = \sum_{h_2} |h_2| |\gamma_{h_1,h_2} |$, $S_{\cdot,
h_2} = \sum_{h_1} |h_1| |\gamma_{h_1,h_2} |$,\vadjust{\goodbreak} and $S_{\cdot, \cdot}
= \sum_{h_1,h_2} |h_1| |h_2| |\gamma_{h_1,h_2}|$ and define the set
\begin{eqnarray*}
&&\mathcal{F}  = \biggl\{ F \dvtx{[-\pi,\pi]}^2 \tends
\RR^{+}, F(\lambda_1, \lambda_2) = \sum
_{h_1,h_2} \gamma_{h_1,h_2} (F) Z_1^{h_1}
Z_2^{h_2},
 S_{h_1, \cdot} < \infty\\
 &&\hspace*{194pt}\forall h_1, S_{\cdot, h_2} <
\infty\ \forall h_2, S_{\cdot, \cdot} < \infty\biggr\}.
\end{eqnarray*}
Note that this class excludes spectra with zeroes, which is a minor
imposition in practice.

In this paper we take Brillinger's approach to asymptotic derivations,
stipulating summability conditions on higher
cumulants of the spatial field. Let us denote an integer-valued
bivariate index by $t \in\ZZ^2$, which has integer coordinates
$(t_1, t_2)$. Then a collection of spatial variables can be written
$\{ \WW_{t^{(1)}}, \WW_{t^{(2)}}, \ldots\}$. The weak stationarity
condition stipulates that joint moments of such variables only
depend upon differences between indices, $t^{(1)} - t^{(2)} =
(t_1^{(1)} - t_1^{(2)}, t_2^{(1)} - t_2^{(2)})$, etc. If we sum a
function with respect to $t \in\ZZ^2$, the notation refers to a
double sum over $t_1$ and $t_2$. A similar notation is used for
frequencies $\lambda\in{[-\pi, \pi]}^2$, in that $\lambda=
(\lambda_1, \lambda_2)$.

Suppose that spatial data is sampled from a true spatial field with
spectrum $\widetilde{F}$, and that we have a collection of
continuous weighting functions $G_j \dvtx{[-\pi, \pi]}^2 \mapsto
\RR^{+}$. The second cumulant function of the spatial field is the
autocovariance function $\gamma_{h}$ with $h \in\ZZ^2$, whereas the
$(k+1)$th cumulant function is denoted
\[
\gamma_{h^{(1)}, h^{(2)}, \ldots, h^{(k)}} = \operatorname{cum} [ \WW_{t},
\WW_{t+h^{(1)}},
\WW_{t+h^{(2)}}, \ldots, \WW_{t+ h^{(k)}} ].
\]
We require absolute summability of these second and fourth cumulant
functions. Then the fourth-order spectrum is well defined via
\[
\widetilde{FF} \bigl( \lambda^{(1)}, \lambda^{(2)},
\lambda^{(3)} \bigr) = \sum_{h^{(1)},h^{(2)},h^{(3)}}
\gamma_{h^{(1)},h^{(2)},h^{(3)}} e^{
\{ -i \lambda^{(1)} \cdot
h^{(1)} - i \lambda^{(2)} \cdot h^{(2)} - i \lambda^{(3)} \cdot
h^{(3)} \} },
\]
with $\cdot$ denoting the dot product of bivariate vectors. More
regularity can be imposed via the condition
%
\begin{equation}
\label{eq:Brillinger} \sum_{h^{(1)}, h^{(2)}, \ldots, h^{(k)} }
\bigl(1 + \overline
{\bigl|h^{(1)}\bigr|} \overline{\bigl|h^{(2)}\bigr|}
\cdots\overline{\bigl|h^{(k)}\bigr|}
\bigr) \bigl| \gamma_{h^{(1)}, h^{(2)}, \ldots
, h^{(k)}} \bigr| < \infty,
\end{equation}
where $\overline{t}$ denotes the product of the components of $t$.
This will be referred to as Condition $B_k$, for any $k \geq1$;
note that $B_2$ implies the summability conditions of the set
$\mathcal{F}$. Finally, recall that the periodogram is computed
from a sample of size $N^2 = N_1 \cdot N_2$. When the regressors
are correctly specified, we will write $\widetilde{\beta}$ for the
true parameter. Then $I_{\widetilde{\beta}} $ denotes the
periodogram of the data $Y$ correctly adjusted for mean effects;
equivalently, it is the periodogram of $\WW^N$.

In addition to assuming that the regressors are correctly specified,
with $\widetilde{\beta}$ the true regression parameter and
$\widetilde{X}$ the regression matrix, we require the following key
assumptions.

\begin{ass*}
\begin{longlist}[(A3)]
\item[(A1)] $\widetilde{F} \in\mathcal{F}$.

\item[(A2)] The spectral density $F_{\theta}$ is twice
continuously differentiable and uniformly bounded above and away
from zero, and moreover all components of $F_{\theta}$, $\nabla
F_{\theta}$, $\nabla\nabla^{\prime} F_{\theta}$ are in
$\mathcal{F}$.

\item[(A3)] The process is weakly stationary of order $k$,
and the Brillinger conditions $B_k$ (\ref{eq:Brillinger}) hold for
all $k \geq1$.

\item[(A4)] The pseudo-true value $\widetilde{\theta}$
exists uniquely in the interior of the parameter space.

\item[(A5)] $H(\theta) = \nabla\nabla^{\prime} \operatorname{KL}
(F_{\theta}, \widetilde{F}) $ is invertible at
$\widetilde{\theta}$.
\end{longlist}
\end{ass*}

Conditions ({A1}), ({A3}) and ({A5}) cannot be
verified from
data, but some assumptions of this nature must be made to obtain
asymptotic formulas. Condition ({A2}) will hold for cepstral
models (and other random field models as well) by the following
argument. The coefficients of the causal and skew fields will have
exponential decay in either index argument, by extensions of the
classical time series argument (see, e.g.,
Hurvich \cite{Hurv:mult:2002}) applied to (\ref{eq:betaRecurs}) and
(\ref{eq:rhoRecurs}). [The time series argument can be directly
applied to (\ref{eq:deltaRecurs}) and (\ref{eq:etaRecurs}) as well.]
Combining these results using (\ref{eq:acfExpress}), the acf of the
cepstral field will also have exponential decay so that $F_{\theta}
\in\mathcal{F}$. Of course, another way to verify this condition is
to examine the boundedness of partial derivatives of the spectrum;
at once we see that ({A2}) holds for the cepstral model, as it
does for moving average random fields.

Although condition ({A4}) may be problematic for certain moving
average models (which may have complicated constraints on
coefficients), the cepstral model uses no constraints on $\theta$,
because the distinct entries of $\Theta$ can be any real number,
independently of all other distinct entries. Euclidean space is open,
so any
pseudo-true value is necessarily contained in the interior. Also,
existence of a pseudo-true value is guaranteed by convexity of the
KL discrepancy.

For the result on Bayesian estimation, we will assume that the model
is correctly specified; the model must also be identifiable, that is,
$F_{\theta_1} = F_{\theta_2}$ implies $\theta_1 = \theta_2$, which
helps ensure asymptotic concentration of the likelihood. We assume
the parameters belong to some compact subset of Euclidean space, and
the true parameter vector lies in the interior. This assumption
can often be accomplished by prior transformation (and is easily
accomplished for the cepstral coefficients in the cepstral model).
Also define the matrix 2-norm of a matrix $A$ via the notation ${\|
A \| }_2$.

\subsection{Technical results}

We begin with an important lemma that extends Lemma 4.1.2 of
Taniguchi and
Kakizawa \cite{taniguchi2000asymptotic} to the spatial context.
%
%
\begin{Lemma}
\label{lem:toeplitz} Let $\Sigma(F_j)$ and $\Sigma(G_j)$ be
block-Toeplitz matrices with $F_j,\break  G_j^{-1} \in\mathcal{F}$ for $1
\leq j \leq m$. Assuming that $N_* = \min\{ N_1, N_2 \} \tends
\infty$, and $N = \sqrt{N_1 \cdot N_2}$,
\[
N^{-2} \operatorname{tr} \Biggl\{ \prod_{j=1 }^m
\Sigma(F_j) \Sigma^{-1} (G_j) \Biggr\} = \Biggl
\langle\prod_{j=1}^m F_j G_j^{-1}
\Biggr\rangle+ O\bigl(N^{-2}\bigr).
\]
\end{Lemma}
Next, we discuss a lemma that provides a central limit theorem for
weighted averages of the spatial periodogram, which is a natural
extension of Lemma 3.1.1 of Taniguchi and
Kakizawa \cite{taniguchi2000asymptotic}. Define
the bias-correction quantities
%
\begin{equation}
\label{eq:biasCorrect} B_1 (\lambda) = \sum_{h_1, h_2}
|h_1| \gamma_{h_1, h_2} e^{-i
\lambda\cdot h},\qquad B_2 (
\lambda) = \sum_{h_1, h_2} |h_2|
\gamma_{h_1, h_2} e^{-i
\lambda\cdot h},
\end{equation}
and use $\langle\!\langle g (\lambda^{(1)}, \lambda^{(2)})
\rangle\!
\rangle$ as a short hand for ${(2 \pi)}^{-4} \int_{{[-\pi,\pi]}^4}
g (\lambda^{(1)}, \lambda^{(2)}) \,d\lambda^{(1)} \,d\lambda^{(2)}$.
%

\begin{Lemma}
\label{lem:Fourier} Assume that $N_* = \min\{ N_1, N_2 \} \tends
\infty$ and let $N = \sqrt{N_1 \cdot N_2}$. Suppose assumption
\textup{({A3})} holds, and that $G_j$ for $1 \leq j \leq J$ are continuous
functions. Let $G_j^* (\lambda) = G_j (-\lambda)$. Then:
\begin{longlist}[(ii)]
\item[(i)] For the unbiased acf estimators, as $N_* \tends\infty$, $
\langle G_j \widehat{I}_{\widetilde{\beta}} \rangle- { \langle G_j
\widehat{I}_{\widetilde{\beta}} \rangle}_N \convinp0$ and $\langle
G_j \widehat{I}_{\widetilde{\beta}} \rangle\convinp\langle G_j
\widetilde{F} \rangle$ for any $1 \leq j \leq J$. Also
\[
N { \bigl\{ \bigl\langle G_j (\widehat{I}_{\widetilde{\beta}} -
\widetilde{F} ) \bigr\rangle\bigr\} }_{j=1}^J \,\convinlaw\,
\mathcal{N} (0, V),
\]
where the covariance matrix $V$ has $jk$th entry
\[
\bigl\langle\bigl\langle G_j G_k \widetilde{FF} \bigl(
\lambda^{(1)}, - \lambda^{(2)}, \lambda^{(2)}\bigr)
\bigr\rangle\bigr\rangle+ \bigl\langle\bigl( G_j G_k^*
+ G_j G_k \bigr) \widetilde{F}^2 \bigr
\rangle.
\]
\item[(ii)] For the biased acf estimators, $\langle G_j I_{\widetilde{\beta}}
\rangle- { \langle G_j I_{\widetilde{\beta}}
\rangle}_N \convinp0$ and $\langle G_j I_{\widetilde{\beta}}
\rangle
\convinp\langle G_j \widetilde{F}
\rangle$ for any $1 \leq j \leq J$. Also, for the same $V$ given
in case \textup{(i)},
\[
N { \bigl\{ \bigl\langle G_j \bigl( I_{\widetilde{\beta}} -
\widetilde{F} + N^{-1} B_1 + N^{-1}
B_2 \bigr) \bigr\rangle\bigr\} }_{j=1}^J
\,\convinlaw\,\mathcal{N} (0,V),
\]
where the bias correction terms are defined in
(\ref{eq:biasCorrect}).
\end{longlist}
\end{Lemma}

The last assertion of Lemma \ref{lem:Fourier} means that utilizing
$I_{\widetilde{\beta}}$ instead of $\widehat{I}_{\widetilde{\beta}}$
will require a bias correction; cf. Guyon \cite{guyon1982parameter}. Both
lemmas are important preliminary results for our main theorems, but
also are of interest in their own right, extending known time
series
results to the spatial context. Although generalizations to
dimensions higher than two seem feasible, the actual mechanics
become considerably more technical. We now state the limit theorems
for our parameter estimates. For the QMLE estimates, we suppose
that they are either exact or approximate Whittle estimates defined
using the unbiased acf estimates.

%
\begin{thmm}
\label{thmm:AsympFreq}
Assume that conditions \textup{({A1})--({A5})} hold and that the
regressors are correctly specified with ${( \widetilde{X}^{\prime}
\widetilde{X} )}^{-1} \tends0$ as $N_* = \min\{ N_1, N_2 \} \tends
\infty$. Then in the case of MLE or the QMLE, both $\widehat{\beta}$
and $\widehat{\theta}$ are jointly asymptotically normal with
distributions given by
\begin{eqnarray*}
N ( \widehat{\theta} - \widetilde{\theta} )& \convinlaw& \mathcal
{N} \bigl( 0,
H^{-1} (\widetilde{\theta}) V (\widetilde{\theta}) H^{\dag} (
\widetilde{\theta}) \bigr),
\\
H(\theta) & =& \nabla\nabla^{\prime} \operatorname{KL} (F_{\theta},
\widetilde{F}),
\\
V(\theta) & =& 2 \bigl\langle\widetilde{F}^2 \nabla
F_{\theta}^{-1} \nabla^{\prime} F_{\theta}^{-1}
\bigr\rangle+ \bigl\langle\bigl\langle\nabla F^{-1}_{\theta}
\nabla^{\prime} F^{-1}_{\theta} \widetilde{FF} \bigl(
\lambda^{(1)}, - \lambda^{(2)}, \lambda^{(2)} \bigr)
\bigr\rangle\bigr\rangle,
\end{eqnarray*}
where $\dag$ denotes an inverse transpose and $N = \sqrt{N_1 \cdot
N_2}$. Also $N ( \widehat{\beta} - \widetilde{\beta} )$ is
asymptotically normal with mean zero and covariance matrix
\[
M^{-1}_X (\widetilde{\theta}) { \bigl[
\widetilde{X}^{\prime} \Sigma\bigl(F^{-1}_{\widetilde{\theta
}}\bigr)
\Sigma(\widetilde{F}) \Sigma\bigl(F^{-1}_{\widetilde{\theta
}}\bigr)
\widetilde{X} \bigr] } M^{\dag}_X (\widetilde{\theta}),
\]
where $M_X (\theta) = \widetilde{X}^{\prime}
\Sigma(F_{ \theta}^{-1}) \widetilde{X}$. Finally, $\widehat{\beta}$
and $\widehat{\theta}$ are asymptotically
independent if the third cumulants of the process are zero.
\end{thmm}
%
%
\begin{Remark}
For a Gaussian process, third and fourth cumulants are zero,
which implies that regression and model parameter estimates are
asymptotically independent, and that $V$ has a simpler form, being
given just by $2 \langle\widetilde{F}^2 \nabla F_{\theta}^{-1}
\nabla^{\prime} F_{\theta}^{-1} \rangle$.
\end{Remark}

%
\begin{Remark}
Application of the same techniques in the case of a one-dimensional
random field, or time series, yields asymptotic normality of
regression and time series parameters under Brillinger's conditions.
To our knowledge, the only other results of this flavor for time
series with regression effects is the work of Pierce \cite{pierce1971},
which focuses on ARIMA models but allows for skewed non-Gaussian
distributions.
\end{Remark}

%
\begin{thmm}
\label{thmm:Bayes} Assume that the data process is Gaussian and
\textup{({A2})} and \textup{({A4})} hold, and that the model is correctly
specified and
is identifiable. Also suppose that the regressors are correctly
specified, with $N^{-2} \widetilde{X}^{\prime} \Sigma^{-1}
(F_{\theta}) \widetilde{X} \tends M(\theta)$ for some $M(\theta)$
satisfying $0 < \sup_{\theta} { \| M( \theta) \| }_2 < \infty$. Then
$\widehat{\beta}_B \convinp\widetilde{\beta}$ and
$\widehat{\theta}_B \convinp\widetilde{\theta}$ as $N_* = \min\{
N_1, N_2 \} \tends\infty$.
\end{thmm}

It is worth comparing the conditions of the two theorems. In
Theorem \ref{thmm:Bayes} the assumption of a correct model makes
({A1}) automatic, and the Gaussian assumption makes ({A3})\vadjust{\goodbreak}
automatic.
Furthermore, the assumption in Theorem \ref{thmm:Bayes} on the
parameters---together with the assumption of a correct
model---automatically entails ({A4}) as well. Theorem \ref
{thmm:AsympFreq}
also assumes ({A5}), which is chiefly needed to establish
asymptotic normality of the frequentist estimates. The Bayesian
result requires a slightly stronger assumption on the regression
matrix in order to get asymptotic concentration of the likelihood.
For example, if we seek to estimate a constant mean by taking
$\widetilde{X}$ to be a column vector of all ones, then $M(\theta)$
exists and is just the scalar $F_{-\theta} (0,0)$; this will be
bounded away from zero and infinity in the cepstral model if all the
cepstral coefficients are restricted to a range of values.

\section{Simulation study}\label{sec:sim}
To demonstrate the effectiveness of our approach, we conducted a
small simulation study using maximum likelihood estimation as
outlined in Sections~\ref{sec:model} and \ref{sec:fitting}. The
model autocovariances were calculated according~(\ref{eq:acfApprox}), with $M=1000$ and $p_1=p_2=2$. The exact
parameter values for the simulation were calibrated to the straw
yield data analysis presented in the Supplement's Appendix A. Grid
sizes of $(15\times15)$, $(20\times20)$, $(20\times25)$ and
$(25\times25)$ were considered, where $(20\times25)$ constitutes
the size grid in our real-data example.

\begin{sidewaystable}
\tablewidth=\textwidth
\caption{Simulation results for the simulation
presented in Section~\protect\ref{sec:sim} ($p_1=p_2=2$). Note, there were
200 simulated datasets and $\overline{\sigma}_\theta$ denotes the
mean standard deviation for parameters $\theta_1,\ldots,\theta_{12}$
for a given simulation (over the 200 datasets). Recall that, for
$j=1,\ldots,12$, the asymptotic standard error for $\theta_j$ equals
$1/N$. Therefore, $\overline{\sigma}_\theta-1/N$ represents the average
difference between the estimated and asymptotic standard error for
$\theta_1,\ldots,\theta_{12}$. The values in the table below are only
reported to three decimal places and the
elements of $\theta$ are described in Section \protect\ref{sec:sim}.
Note that $N=\sqrt{N_1\cdot N_2}$ and that $\theta_{13}=\Theta_{0,0}$}
\label{tab:sim}
\begin{tabular*}{\textwidth}{@{\extracolsep{\fill}}ld{2.3}d{2.3}ccd{2.3}ccd{2.3}ccd{2.3}cc@{}}
\hline
&&\multicolumn{12}{c@{}}{$\bolds{N_1\times N_2}$}\\[-4pt]
&&\multicolumn{12}{c@{}}{\hrulefill}\\
&&\multicolumn{3}{c}{$\bolds{15\times15}$}&\multicolumn{3}{c}{$\bolds{20\times
20}$}&\multicolumn{3}{c}{$\bolds{20 \times25}$}& \multicolumn{3}{c@{}}{$\bolds{25
\times25}$}\\[-4pt]
&&\multicolumn{3}{c}{\hrulefill}&\multicolumn{3}{c}{\hrulefill}&
\multicolumn{3}{c}{\hrulefill}& \multicolumn{3}{c@{}}{\hrulefill}\\
&\multicolumn{1}{c}{\textbf{True}} & \multicolumn{1}{c}{\textbf{Mean}}& \textbf{SD} &
 \textbf{MSE} &
\multicolumn{1}{c}{\textbf{Mean}}& \textbf{SD} & \textbf{MSE} &
\multicolumn{1}{c}{\textbf{Mean}}& \textbf{SD} & \textbf{MSE} &
\multicolumn{1}{c}{\textbf{Mean}} & \textbf{SD}&
\multicolumn{1}{c@{}}{\textbf{MSE}}\\
\hline
$\theta_1$ & 0.009 &-0.019& 0.084& 0.008 & 0.000& 0.056&0.003 &-0.001&
0.049 &0.003 &-0.002& 0.041& 0.002 \\
$\theta_2$ &-0.028 &-0.037& 0.079& 0.006 &-0.029& 0.057&0.003 &-0.033&
0.048 &0.002 &-0.039& 0.044& 0.002 \\
$\theta_3$ & 0.132 & 0.123& 0.081& 0.007 & 0.127& 0.059&0.003 & 0.133&
0.044 &0.002 & 0.129& 0.045& 0.002 \\
$\theta_4$ & 0.067 & 0.054& 0.080& 0.007 & 0.057& 0.059&0.004 & 0.058&
0.047 &0.002 & 0.063& 0.049& 0.002 \\
$\theta_5$ & 0.271 & 0.266& 0.078& 0.006 & 0.265& 0.056&0.003 & 0.269&
0.052 &0.003 & 0.259& 0.043& 0.002 \\
$\theta_6$ & 0.383 & 0.367& 0.074& 0.006 & 0.370& 0.057&0.003 & 0.377&
0.047 &0.002 & 0.379& 0.038& 0.001 \\
$\theta_7$ & 0.001 &-0.006& 0.079& 0.006 &-0.001& 0.060&0.004 &-0.005&
0.050 &0.003 &-0.001& 0.044& 0.002 \\
$\theta_8$ &-0.017 &-0.028& 0.077& 0.006 &-0.023& 0.059&0.004 &-0.022&
0.049 &0.002 &-0.020& 0.045& 0.002 \\
$\theta_9$ &-0.003 &-0.023& 0.082& 0.007 &-0.012& 0.052&0.003 &-0.009&
0.047 &0.002 &-0.009& 0.046& 0.002 \\
$\theta_{10}$ &-0.055 &-0.090& 0.079& 0.008 &-0.064& 0.049&0.002
&-0.053& 0.053 &0.003 &-0.061& 0.040& 0.002 \\
$\theta_{11}$ &-0.015 &-0.035& 0.085& 0.008 &-0.022& 0.054&0.003
&-0.021& 0.048 &0.002 &-0.021& 0.044 & 0.002 \\
$\theta_{12}$ & 0.144 & 0.129& 0.072& 0.005 & 0.134& 0.053&0.003 &
0.138& 0.047 &0.002 & 0.137& 0.042 & 0.002 \\
$\theta_{13}$ &-0.871 &-0.968& 0.096& 0.019 &-0.922& 0.080&0.009
&-0.902& 0.064 &0.005 &-0.902& 0.060& 0.004 \\
$\beta_0$ & 7.646 & 7.632& 0.275& 0.076 & 7.651& 0.200&0.040 & 7.627&
0.167 &0.028 & 7.645& 0.173& 0.030 \\
$\beta_1$ &-0.035 &-0.035& 0.020& 0.000 &-0.034& 0.011&0.000 &-0.034&
0.011 &0.000 &-0.034& 0.008& 0.000 \\
$\beta_2$ &-0.059 &-0.056& 0.024& 0.001 &-0.060&
0.013&0.000 &-0.059& 0.009&0.000 &-0.060& 0.009& 0.000
\\[3pt]
\multicolumn{2}{@{}l}{$\overline{\sigma}_\theta
-(1/N)$}&\multicolumn
{3}{c}{0.0124}&\multicolumn{3}{c}{0.0059}&\multicolumn
{3}{c}{0.0037}&\multicolumn{3}{c}{0.0033}\\
\hline
\end{tabular*}
\end{sidewaystable}

For this simulation, we generated 200 Gaussian datasets with
parameters $\theta= J \Theta$ corresponding to quadrants I and II
of the grid $\Theta$,\setcounter{footnote}{1}\footnote{That is, $\theta= [ \Theta_{-2,2},
\Theta_{-2,1}, \Theta_{-1,2}, \Theta_{-1,1}, \Theta_{0,2},
\Theta_{0,1}, \Theta_{1,2}, \Theta_{1,1}, \Theta_{2,2},
\Theta_{2,1}, \Theta_{-2,0},  \Theta_{-1,0},\break  \Theta_{0,0}
]^{\prime} =(\theta_1,\theta_2,\ldots,\theta_{13})^{\prime}$.}
and $\beta=(\beta_0,\beta_1,\beta_2)'$; see Table~\ref{tab:sim}. In
this case, $\beta_1$ and $\beta_2$ correspond to ``row'' and
``column'' effects, respectively, in the agricultural experiment
considered. Here, the row and column effects are obtained by
regressing the vectorized response on the corresponding row and
column indices (since rows and columns are equally spaced). The
$\widetilde{X}$ matrix used in this simulation consisted of a column
of ones followed by columns associated with the row and column
effects and was taken from the analysis presented in the
Supplement's Appendix A. Given $\theta$ and $\beta$, we simulate
directly from the corresponding multivariate Gaussian distribution.
However, in cases where the grid size is extremely large, another
potential approach to simulation would be circular embedding
(Chan and Wood \cite{chan1999simulation}, Wood and Chan \cite
{wood1994simulation}), though
it would be
necessary to properly account for any regression effects. The log
Gaussian likelihood (up to constants) given by (\ref{eq:GaussLik})
was numerically maximized for each simulated dataset using the
\textit{optim} function in R (R Development Core Team
\cite{R:software}).\looseness=-1

As demonstrated in Table~\ref{tab:sim}, through an assessment of
mean square error (mse), the model parameters can be estimated with
a high degree of precision. Additionally, Table~\ref{tab:sim}
illustrates that our asymptotic theory agrees with the finite sample
estimates for different grid sizes. Specifically, we provide the
difference between the mean standard deviation (over all of the
cepstral parameters, except $\Theta_{0,0}$)\vadjust{\goodbreak} and the asymptotic
standard deviation. This simulation shows that, as the grid size
increases, the difference between the estimated standard error and
the asymptotic standard error goes to zero on average. We also
provide the mean, standard deviation and mse for the individual
parameters, including the mean parameters $\beta$; this demonstrates
the bias properties, as well as the fact that the mse goes to zero
as the grid size increases. Finally, the average $p$-value for the
Shapiro--Wilks test of normality for each simulation grid size (over
all of the cepstral parameters) was greater than 0.4, with only one
parameter out of the thirteen cepstral parameters from each simulation
not exhibiting
normality. Hence, the estimated parameters converge to their
asymptotic distribution and, as expected, their precision increases
with sample size.

\section{Conclusion}\label{sec:conclusion}
The general modeling approach and asymptotic theory we propose
extends the spatial random field literature in several directions.
By providing recursive formulas for calculating autocovariances,
from a given cepstral random field model, we have facilitated usage
of these models in both Bayesian and likelihood settings. This
is extremely notable as many models suffer from a constrained
parameter space, whereas the cepstral random field model imposes no
constraints on the parameter values. More specifically, the
autocovariance matrix obtained from our approach is guaranteed to be
positive definite.

In addition, we establish results on consistency and asymptotic
normality for an expanding domain. This provides a rigorous
platform for conducting model selection and statistical inference.
The asymptotic results are proven generally and can be viewed as an
independent contribution to the random field literature, expanding
on the results of Mardia and Marshall \cite{mardia1984maximum} and
others, such as
Guyon \cite{guyon1982parameter}.
The simulation results support the theory, and the methods are
illustrated through an application to straw yield data from an
agricultural field experiment (Supplement's Appendix A). In this
setting, it is readily seen
that our model is easily able to characterize the underlying spatial
dependence structure.

\section*{Acknowledgments}
We thank an Associate Editor and two anonymous referees for providing detailed
comments that helped substantially improve this article. The authors
would also like to thank Aaron Porter for his assistance with
implementation of the Moran's I statistic. This article is released to
inform interested parties of research and to encourage discussion. The
views expressed on statistical issues are those of the authors and not
necessarily those of the U.S. Census Bureau. See \cite{supp} for the
supplementary material.


\begin{supplement}[id=suppA]
\stitle{Supplement to asymptotic theory of cepstral random fields\\}
\slink[doi]{10.1214/13-AOS1180SUPP} 
\sdatatype{.pdf}
\sfilename{aos1180\_supp.pdf}
\sdescription{The supplement contains a description of further
applications of the cepstral model, analysis of straw yield data, as
well as all proofs.}
\end{supplement}

%

%

\printaddresses

\end{document}